\documentclass[11pt,a4paper]{article}
\usepackage[utf8]{inputenc}
\usepackage[english]{babel}
\usepackage[T1]{fontenc}
\usepackage{lmodern}
\usepackage{amsmath}
\usepackage{accents}
\usepackage{amsfonts}
\usepackage{amssymb}
\usepackage{amsthm}
\usepackage{amscd}
\usepackage{url}
\usepackage{enumitem}
\usepackage{hhline}
\usepackage{dsfont}
\usepackage{array}
\usepackage{geometry}
\usepackage{tikz}
\usetikzlibrary{decorations.pathreplacing}
\usetikzlibrary{patterns}
\usepackage{xcolor}
\usetikzlibrary{arrows}
\usepackage{mdframed}
\usepackage{caption}
\usepackage[titletoc]{appendix}
\usepackage{csquotes}
\usepackage{hyperref}
\usepackage{float}
\usepackage{bm}
\usepackage{microtype}
\usepackage{graphicx}

\mdfsetup{nobreak=true}

\AtBeginDocument{
  \label{CorrectFirstPageLabel}
  
}

\usepackage[backend=biber]{biblatex}
\newcolumntype{L}{>{$}l<{$}}

\newcommand\N{\mathbb{N}}
\newcommand\Nn{\mathcal{N}}

\newcommand\Z{\mathbb{Z}}

\newcommand\R{\mathbb{R}}
\newcommand\C{\mathbb{C}}

\newcommand\E{\mathbb{E}}

\DeclareMathOperator\Var{Var}
\newcommand\dd{\mathrm{d}}
\newcommand\CC{\mathcal{C}}

\DeclareMathOperator\Cov{Cov}
\newcommand\Card{\mathrm{Card}}

\newcommand\T{\mathbb{T}}

\newcommand{\enstq}[2]{\left\{#1~\middle|~#2\right\}}

\DeclareMathOperator\dist{dist}
\DeclareMathOperator\sinc{sinc}

\DeclareMathOperator\diag{diag}
\DeclareMathOperator\off{off}
\DeclareMathOperator\Law{Law}
\newcommand\one{\mathds{1}}

\newcommand\numberthis{\addtocounter{equation}{1}\tag{\theequation}}
\newcommand\jump{\par\medskip}
\newcommand\quand{\quad\text{and}\quad}

\theoremstyle{plain}
\newtheorem{theo}{Theorem}[section]
\newenvironment{theorem}%
  {\begin{mdframed}[backgroundcolor=white]\begin{theo}}%
  {\end{theo}\par\vspace{0.1cm}\end{mdframed}}

\theoremstyle{plain}
\newtheorem{prop}[theo]{Proposition}
\newenvironment{proposition}%
  {\begin{mdframed}[backgroundcolor=white]\begin{prop}}%
  {\end{prop}\par\vspace{0.1cm}\end{mdframed}}

\theoremstyle{plain}
\newtheorem{coro}[theo]{Corollary}
\newenvironment{corollary}%
  {\begin{mdframed}[backgroundcolor=white]\begin{coro}}%
  {\end{coro}\par\smallskip\end{mdframed}}
  
\theoremstyle{plain}
\newtheorem{lemm}[theo]{Lemma}
\newenvironment{lemma}%
  {\begin{mdframed}[backgroundcolor=white]\begin{lemm}}%
  {\end{lemm}\par\vspace{0.cm}\end{mdframed}}

\theoremstyle{definition}
\newtheorem{remark}[theo]{Remark}

\makeatletter
\def\blfootnote{\gdef\@thefnmark{}\@footnotetext}
\makeatother
  
\makeatletter
\renewcommand*\env@matrix[1][*\c@MaxMatrixCols c]{%
  \hskip -\arraycolsep
  \let\@ifnextchar\new@ifnextchar
  \array{#1}}
\makeatother

\begin{document}

\title{\Huge{Variance of the number of zeros of dependent Gaussian trigonometric polynomials}}
\author{Louis Gass}
\maketitle
\blfootnote{Univ Rennes, CNRS, IRMAR - UMR 6625, F-35000 Rennes, France.}
\blfootnote{This work was supported by the ANR grant UNIRANDOM, ANR-17-CE40-0008.}
\blfootnote{Email: louis.gass@ens-rennes.fr}

\begin{abstract} We compute the variance asymptotics for the number of real zeros of trigonometric polynomials with random dependent Gaussian coefficients and show that under mild conditions, the asymptotic behavior is the same as in the independent framework. In fact our proof goes beyond this framework and makes explicit the variance asymptotics of various models of random Gaussian processes. Our proof relies on intrinsic properties of the Kac--Rice density in order to give a short and concise proof of variance asymptotics.
\end{abstract}

\setcounter{tocdepth}{2}
\tableofcontents\jump

\newpage
\section{Introduction}
The asymptotic behavior of the variance of the number of zeros of random trigonometric polynomials $\sum a_k\cos(kt) + b_k\sin(kt)$ with independent Gaussian coefficients has been established in \cite{Gra11}. Since then, the variances of numerous models have been studied: for instance, see \cite{Aza16} for the analogous model $\sum a_k\cos(kt)$, or more recently \cite{Lub21} for random orthogonal polynomials on the real line. In this paper we make explicit the asymptotics of the variance of the number of zeros of random trigonometric polynomials with dependent coefficients.\jump

This paper serves three purposes. Firsty, it is the natural continuation of \cite{Ang19} and \cite{Ang21b}, in which the asymptotic behavior of the expected number of zeros of this model has been established. Secondly, we avoid the use of the explicit expression for the first and second order Kac densities for the number of zeros (see e.g. \cite[Lemma 2.2]{Lub21}), but we rather exploit on one hand the intrinsic properties of the Kac density and on the other hand the convergence in distribution of the model of random trigonometric polynomials towards the $\sinc$ process. This point of view allows us to avoid heavy computations that usually comes with the Kac--Rice formula, in particular in the non-stationary case, and allows us to give a unified point of view for all the previously cited Gaussian models for which the variance asymptotics is known. \jump


Let us now detail our model of random trigonometric polynomial with dependent coefficients. Let $\T = \R/2\pi\Z$ be the one-dimensional torus, which can be identified with a segment of $\R$ of length $2\pi$. For $s,t\in \T$ we define the distance $\dist(s,t)=d(s-t,2\pi\Z)$. For $s\in\T$ we define 
\[X_n(s) = \frac{1}{\sqrt{n}}\sum_{k=1}^{n}a_k\cos(ks) + b_k\sin(ks),\numberthis\label{eqbis:17}\]
with $(a_k)_k, (b_k)_k$ two independent stationary centered Gaussian processes with correlation function $\rho : \Z\rightarrow \R$. That is for $k,l\geq 0$, $\E[a_ka_l] = \E[b_kb_l] = \rho(k-l)$. Thanks to Bochner Theorem, the correlation function $\rho$ is associated with a spectral measure $\mu$ on the torus $\T$ via the relation
\[\rho(k) = \frac{1}{2\pi}\int_\T e^{-iku}\dd \mu(u).\]
We denote by $Z_{X_n}(I)$ the number of zeros of $X_n$ on a subinterval $I$ of the torus $\T$. Under suitable conditions on the spectral measure $\mu$, it has been shown in \cite{Ang19} and \cite{Ang21} that the expectation of the number of zeros of the process $X_n$ on $\T$ behaves like $\frac{2}{\sqrt{3}}n$, as in the independent framework. The following theorem makes explicit the variance asymptotics for the number of zeros of the process $X_n$ on $\T$, as $n$ grows to infinity.
\begin{theorem}
\label{theorembis1}
We suppose that the spectral measure $\mu$ has a positive continuous density $\psi$ with respect to the Lebesgue measure on the torus $\T$. Then there is an explicit positive constant $\gamma_2$ that does not depend on $\psi$, such that for any subinterval $I$ of the torus,
\[\lim_{n\rightarrow +\infty} \frac{\Var(Z_{X_n}(I))}{n} = \mathrm{length}(I)\,\gamma_2.\] 
\end{theorem}

The universal constant $\gamma_2 \simeq 0.089$ is thus the same constant computed in \cite{Gra11} in the particular framework of independent Gaussian random variables. Note that this universality result is not unsurprising, given the non-universality behavior in the non Gaussian framework (see \cite{Bal19}). \jump

Theorem \ref{theorembis1} is in fact a corollary of the next Theorem \ref{theorembis2}, that puts forwards the principal ingredients necessary to obtain such a universal asymptotics for the variance. Let $(X_n)_n$ be a sequence of centered Gaussian random processes defined on an open subinterval $I$ of $\T$ or $\R$. We define for $n\in\N$ the process $Y_n = X_n(\,.\,/n)$, defined on the open subinterval $nI$ of the rescaled torus $n\T$ or $\R$ endowed with its canonical distance. We denote $r_n$ the covariance function of $Y_n$
\[\forall s,t\in nI,\qquad r_n(s,t) := \E[Y_n(s)Y_n(t)].\]
Let also $Y_\infty$ be a stationary Gaussian random process on $\R$ with covariance function $r_\infty$, and $\psi$ be a uniformly continuous function on $I$, bounded above and below by positive constants.\jump

We denote by $C_0(\R)$ the space of continuous functions that converges to zero in $\pm\infty$.
\begin{theorem}
\label{theorembis2}
We suppose that the sequence of random processes $(X_n)_n$ satisfies the following two conditions.
\begin{itemize}
\item[$(H_1)$] The sequence of processes $(X_n)_{n\in\N}$ is of class $\CC^2$ on $I$ and for $u,v\in\{0,1,2\}$, the following convergence holds uniformly for $x\in I$ and locally uniformly for $s,t\in \R$,
\[\lim_{n\rightarrow+\infty} r_n^{(u,v)}\left(nx+s,nx+t\right) = \psi(x)r_\infty^{(u,v)}(s,t).\]
\item[$(H_2)$] There is a positive function $g\in L^2(\R)\cap C_0(\R)$ such that for $u,v\in\{0,1\}$,
\[\forall s,t\in nI,\quad|r_n^{(u,v)}(s,t)|\leq g(\dist(s,t)).\]
\end{itemize}
Then there is an explicit positive constant $\gamma_2$ depending only on $r_\infty$, such that
\[\lim_{n\rightarrow +\infty} \frac{\Var(Z_{X_n}(I))}{n} = \mathrm{length}(I)\,\gamma_2.\] 
\end{theorem}

This last Theorem \ref{theorembis2} covers Theorem \ref{theorembis1}, as we will prove that the assumptions on the spectral measure $\mu$ in Theorem \ref{theorembis1} ensure that the associated sequence of processes $(X_n)_{n}$ satisfies hypotheses $(H_1)$ and $(H_2)$ of Theorem \ref{theorembis2}. \jump

In fact, Theorem \ref{theorembis2} covers various previously known results about the asymptotics of variance of the number of zeros of general random trigonometric polynomials. For instance, Theorem \ref{theorembis2} allows to make explicit the variance asymptotics for the number of zeros on any compact subinterval of $]0,\pi[$ of the process 
\[\widetilde{X}_n(s)= \frac{1}{\sqrt{n}}\sum_{k=0}^n a_k\cos(ks+\theta),\quad\quad\theta\in\R,\numberthis\label{eqbis:33}\]
with $(a_k)_{k\geq 0}$ a stationary sequence of centered random Gaussian variables whose spectral measure has a positive continuous density. Note that the variance asymptotics for the number of zeros on $[0,\pi]$ for this process was established in the case of Gaussian iid $(a_k)_{k\geq 0}$ in \cite{Aza16}. More recently, the authors in \cite{Lub21} established the variance asymptotics for the number of zeros of sums of real orthogonal polynomials with iid Gaussian coefficients. Let $(P_n)_n$ be a sequence of real orthogonal polynomials with respect to a density measure $\dd\phi$ supported on $[-1,1]$. We define
\[X_n(s) = \frac{1}{\sqrt{n}}\sum_{k=0}^n a_kP_k(\cos(s)),\] 
where $(a_n)_{n\geq 0}$ are iid centered Gaussian random variables. Under suitable conditions on the density $\phi$, hypotheses $(H_1)$ and $(H_2)$ of Theorem \ref{theorembis2} are satisfied with $\psi= 1/\phi$ and $Y_\infty$ a process with $\sinc$ covariance function, and Theorem \ref{theorembis2} implies \cite[Cor.~1.3]{Lub21}.\jump

Section \ref{secbis1} of the paper is devoted to the study of the Kac density, that is the integrand in the Kac--Rice formula that gives the second moment of the number of zeros of a general one-dimensional Gaussian process. Under suitable regularity assumptions on the covariance function, we the first remove the apparent singularity of the Kac density near the diagonal. We then show a factorization property of the second order Kac density.\jump

In Section \ref{secbis3}, we make explicit the asymptotics of the Kac density associated with a sequence of processes satisfying the hypotheses of Theorem \ref{theorembis2}. Hypothesis $(H_1)$ implies in particular that as $n$ grows to infinity, the sequence $(Y_n)_n$ has a limit proportional to the stationary process $Y_\infty$ on $\R$. Together with the results of Section \ref{secbis1}, we deduce that the Kac density associated to the process $Y_n$ uniformly converges towards the Kac density of the process $Y_\infty$. This fact leads to the proof of Theorem \ref{theorembis2}. We then check that the model of trigonometric polynomials with dependent coefficients satisfying the hypotheses of Theorem \ref{theorembis1}, also satisfies the hypotheses of Theorem \ref{theorembis2}.

\section{Kac--Rice formula for the second moment}\label{secbis1}
In this section we give a short study of the Kac density for a general centered Gaussian process satisfying some natural regularity assumptions. We first remove the apparent singularity along the diagonal, which is the object of Lemma \ref{lemmabis9} and then prove a factorization property of the Kac density, whose result is contained in Corollary \ref{lemmabis8}.\jump

In the following, let $Y$ be a centered Gaussian process defined on an interval $I$ with covariance function $r$. We assume that the process $Y$ has $\CC^2$ sample paths and that the joint distribution of the Gaussian vector $(Y(s),Y(t))$ is nondegenerate for $s\neq t$. We denote by $p_s$ the density of $Y(s)$ and by $p_{s,t}$ the density of the Gaussian vector $(Y(s),Y(t))$. We will assume the existence of a constant $M$ such that for all $s,t\in I$ and $u,v\leq 2$,
\[|r^{(u,v)}(s,t)|\leq M\quand\det \Cov(Y(s),Y'(s))\geq \frac{1}{M}.\numberthis\label{eqbis:02}\]
\subsection{Kac--Rice formula for the variance}
We define 
\[Z_Y(I) := \Card\enstq{s\in I}{Y(s)=0},\]
the number of zeros of $Y$ in the interval $I$. The following proposition gives an integral representation for the first and second moment of $Z_Y$, and a proof can be found in \cite[Thm.~3.2]{Aza09}.
\begin{proposition}[Kac--Rice formula]
\begin{align*}
\E[Z_Y] &= \int_I \rho_1(s)\dd s \quad\quand\quad \E[Z_Y^2]-\E[Z_Y] = \iint_{I^2} \rho_2(s,t)\dd s\dd t,
\end{align*}
with
\[\rho_1(s) = \E\left[|Y'(s)|\;\middle|\; Y(s)=0\right]p_s(0),\]
\[\rho_2(s,t) = \E\left[|Y'(s)||Y'(t)|\;\middle|\; Y(s)=Y(t)=0\right]p_{s,t}(0,0).\]
\end{proposition}

The function $\rho_1$ (resp. $\rho_2$) is called the Kac density associated to the first (resp. second) moment of the number of zeros. We then have the integral representation
\begin{align*}
\Var(Z_Y) &= \left(\E[Z_Y^2] - \E[Z_Y] - \E[Z_Y]^2\right)  + \E[Z_Y]\\
&=\left(\iint_{I^2} \rho_2(s,t) - \rho_1(s)\rho_1(t)\dd s \dd t\right) + \int_I\rho_1(s)\dd s.\numberthis\label{eqbis:29}
\end{align*}

Note that the quantity $\rho_2(s,s)$ is ill-defined, since the Gaussian vector $(Y(s),Y(s))$ is degenerate. The following lemma allows us to remove the nondegeneracy of the function $\rho_2(s,t)$ when $s$ and $t$ are close, in order to show that the second moment is indeed well-defined. The following procedure is standard, see for instance \cite{Anc20} for a general treatment, or \cite[Prop.~5.8]{Aza09}. The constant $M$ is defined in \eqref{eqbis:02}.
\begin{lemma}
\label{lemmabis9}
There is a positive constant $\eta$ depending only on $M$, such that for $|s-t|\leq \eta$,
\[\rho_2(s,t) \leq M^{3/2}|t-s|.\]
\end{lemma}
\begin{proof}
We define for $s\neq t$ the quantities
\[Y[s,t] = \frac{Y(t)-Y(s)}{t-s}\quand Y[s,s,t] = \frac{Y[s,t]-Y'(s)}{t-s}.\]
And we extend them by continuity respectively by $Y'(s)$ and $Y''(s)/2$ when $s=t$. The mean value theorem and the uniform bounds \eqref{eqbis:02} on $r^{(u,v)}$ imply that
\[\left\|\Cov(Y(s),Y[s,t]) - \Cov(Y(s),Y'(s))\right\|\leq M|t-s|\quand \E[Y[s,s,t]^2]\leq M.\]
The lower bound \eqref{eqbis:02} on $\det \Cov(Y(s),Y'(s))$ implies that for some positive constant $\eta$ depending only on $M$, and $s,t\in I$ such that $|s-t|\leq \eta$,
\[\det\Cov(Y(s),Y[s,t]) \geq \frac{1}{2M}.\]
The density $p_{s,t}(0,0)$ of the vector $(Y(s),Y(t))$ then satisfies for  $|s-t|\leq \eta$
\[p_{s,t}(0,0) = \frac{1}{2\pi\sqrt{\det(\Cov(Y(s),Y(t)))}} = \frac{1}{2\pi|t-s|\sqrt{\det\Cov(Y(s),Y[s,t])}}\leq \frac{\sqrt{M}}{|t-s|}.\]
The conditional Gaussian density appearing in the Kac-Rice formula satisfies, by Cauchy-Schwarz inequality for conditional expectation
\begin{align*}
\E\left[|Y'(s)||Y'(t)|\;\middle|\; Y(s)=Y(t)=0\right]&=|t-s|^2\E\left[|Y[s,s,t]||Y[t,t,s]|\;\middle|\; Y(s)=Y[s,t]=0\right]\\
&\leq |t-s|^2\sqrt{\E\left[Y[s,s,t]^2\right]\E\left[Y[t,t,s]^2\right]}\\
&\leq M|t-s|^2,
\end{align*}
and the conclusion follows.
\end{proof}

\begin{remark}
\label{remark4}
Making explicit the convergence, one has in fact
\[\lim_{t\rightarrow s} \frac{\rho_2(s,t)}{|t-s|} = \frac{\E\left[(Y''(s))^2\;\middle|\; Y(s)=Y'(s)=0\right]}{8\pi\sqrt{\det\Cov(Y(s),Y'(s))}} = \frac{\det(\Cov(Y(s),Y'(s),Y''(s))}{8\pi(\det\Cov(Y(s),Y'(s)))^{3/2}}.\]
\end{remark}
\subsection{Matrix notations}
Before going further, we will need a few matrix notations for the next subsection. In the following, $\Omega$ (resp. $\Sigma$) are square matrices of size $2$ (resp. $4$). We write
\[\Omega = \begin{pmatrix}
\Omega_{11} & \Omega_{12} \\ 
\Omega_{21} & \Omega_{22}
\end{pmatrix} \quand \Sigma = \begin{pmatrix}[c|c]
\Sigma_{11} & \Sigma_{12} \\
\hline
\Sigma_{21} & \vphantom{\int^\int}\Sigma_{22}\quad
\end{pmatrix},
\]
where for $i,j\in\{1,2\}$, $\Omega_{ij}$ are real numbers and $\Sigma_{ij}$ are square matrices of size $2$. We define the diagonal and off diagonal matrices
\[\Omega^{\diag} = \begin{pmatrix}
\Omega_{11} & 0 \\ 
0 & \Omega_{22}
\end{pmatrix},\quad\Omega^{\off} = \begin{pmatrix}
0 & \Omega_{12} \\ 
\Omega_{21} & 0
\end{pmatrix},\quad\Sigma^{\diag} = \begin{pmatrix}[c|c]
\Sigma_{11}^{\diag}\vphantom{\int^\int} & \Sigma_{12}^{\diag} \\
\hline
\Sigma_{21}^{\diag} & \vphantom{\int^\int}\Sigma_{22}^{\diag}\quad
\end{pmatrix},\quad \Sigma^{\off} = \begin{pmatrix}[c|c]
\Sigma_{11}^{\off}\vphantom{\int^\int} & \Sigma_{12}^{\off} \\
\hline
\Sigma_{21}^{\off} & \vphantom{\int^\int}\Sigma_{22}^{\off}\quad
\end{pmatrix},
\]
so that
\[\Omega = \Omega^{\diag} + \Omega^{\off}\quand \Sigma = \Sigma^{\diag} + \Sigma^{\off}.\]
If $\Omega_{11}$ (resp. $\Sigma_{11}$) is non-zero (resp. invertible) we define the Schur complements
\[\Omega_c := \Omega_{22} - \frac{\Omega_{12}\Omega_{21}}{\Omega_{11}}\quand\Sigma_c := \Sigma_{22} - \Sigma_{21}\Sigma_{11}^{-1}\Sigma_{12}.\numberthis\label{eqbis:30}\]
By row reduction, one has
\[\det(\Sigma) = \det(\Sigma_{11})\det(\Sigma_c)\quand (\Sigma^{-1})_{22} = (\Sigma_c)^{-1}\numberthis\label{eqbis:06}.\]
\subsection{Factorization of the Kac density}
In this section, we explicit the Kac densities $\rho_1$ and $\rho_2$ in terms of Gaussian covariance matrices and we prove a useful factorization property for the Kac density $\rho_2$ that will be used in the next section. We define
\[\Omega(s) := \Cov(Y(s),Y'(s))\quand \Sigma(s,t) = \Cov(Y(s),Y(t),Y'(s),Y'(t)).\]
From these definitions and the notation of the previous subsection, one has
\[\det \Sigma^{\diag}(s,t) = \det \Omega(s)\det \Omega(t)\quand (\Sigma^{\diag}(s,t))_c = \begin{pmatrix}
\Omega_c(s) & 0 \\ 
0 & \Omega_c(t)
\end{pmatrix}.\numberthis\label{eqbis:03}\]
The following Lemma \ref{lemmabis7} relies the Schur complement of a matrix with Gaussian conditioning. A proof can be found in \cite{Aza09}.
\begin{lemma}
For $s\neq t$, one has
\label{lemmabis7}
\[\Law((Y'(s))\,|\,Y(s)=0)\sim \Nn(0,\Omega_c(s)).\]
\[\Law((Y'(s),Y'(t))\,|\,Y(s)=Y(t)=0)\sim \Nn(0,\Sigma_c(s,t)).\]
\end{lemma}

We define the function $\tilde{\rho}_1$ (resp. $\tilde{\rho}_2$) on the space of symmetric positive definite matrices of size $2$ (resp. $4$) as
\[\tilde{\rho}_1(\Omega) = \frac{1}{2\pi\sqrt{\det\Omega}}
\int_\R |y|\exp\left(-\frac{y^2}{2\Omega_c}\right)\dd y\]
and 
\[\tilde{\rho}_2(\Sigma) = \frac{1}{(2\pi)^2\sqrt{\det\Sigma }}
\iint_{\R^2} |y_1||y_2|\exp\left(-\frac{1}{2}\,y^T(\Sigma_c)^{-1}y\right)\dd y_1\dd y_2\]
The previous Lemma \ref{lemmabis7} and relations \eqref{eqbis:03} then implies the following formulas
\[\rho_1(s) = \tilde{\rho}_1(\Omega(s)),\quad\rho_2(s,t) = \tilde{\rho}_2(\Sigma(s,t))\quand \rho_1(s)\rho_1(t) = \tilde{\rho}_2(\Sigma^{\diag}(s,t)).\numberthis\label{eqbis:11}\]
When the random Gaussian vectors $(Y(s),Y'(s))$ and $(Y(t),Y'(t))$ are independent then one has directly from Gaussian conditioning the equality $\rho_2(s,t) = \rho_1(s)\rho_1(t)$. The following Lemma \ref{lemmabis4} and Corollary \ref{lemmabis8} makes explicit the error term between these two quantities.
\begin{lemma}
\label{lemmabis4}
Let $K$ be a compact subset of the symmetric positive definite matrices of size $4$. There is a constant $C_K$ depending only on the compact $K$ such that
\[\forall \Sigma\in K,\quad\left|\tilde{\rho}_2(\Sigma) - \tilde{\rho}_2(\Sigma^{\diag})\right|\leq C_K\|\Sigma^{\off}\|^2.\]
\end{lemma}
\begin{proof}
A straightforward computation shows that 
\[\left((\Sigma^{\diag})^{-1}\Sigma^{\off}(\Sigma^{\diag})^{-1}\right)^{\diag}=\left((\Sigma^{\diag})^{-1}\Sigma^{\off}\right)^{\diag} = 0.\numberthis\label{eqbis:04}\]
In particular the trace of these matrices are zero and a Taylor expansion of the determinant function yields
\[\det \Sigma = \det\Sigma^{\diag} + O(\|\Sigma^{\off}\|^2).\numberthis\label{eqbis:05}\]
The expansion of the inverse of a matrix and identity \eqref{eqbis:06} yields
\[(\Sigma_c)^{-1} = \left(\Sigma^{\diag} + \Sigma^{\off}\right)^{-1}_{22} =((\Sigma^{\diag})_c)^{-1} + \left((\Sigma^{\diag})^{-1}\Sigma^{\off}(\Sigma^{\diag})^{-1}\right)_{22} + O(\|\Sigma^{\off}\|^2).\]
By identity \eqref{eqbis:04}, there is an explicit coefficient $h$ and a constant $C_K$ such that
\[\left((\Sigma^{\diag})^{-1}\Sigma^{\off}(\Sigma^{\diag})^{-1}\right)_{22} = \begin{pmatrix}
0 & h \\ 
h & 0
\end{pmatrix}\quad\text{with}\quad |h| \leq C_K\|\Sigma^{\off}\|.\numberthis\label{eqbis:07}\]
We now express the difference
\begin{align*}
\tilde{\rho}_2(\Sigma) - \tilde{\rho}_2(\Sigma^{\diag}) = \frac{1}{4\pi^2}(R_1 + R_2),
\end{align*}
with
\[R_1 = \left[\frac{1}{\sqrt{\det \Sigma }}\!-\!\frac{1}{\sqrt{\det \Sigma^{\diag}}}\right]
\left(\iint_{\R^2} |y_1||y_2|\exp\left(-\frac{1}{2}\,y^T(\Sigma_c)^{-1}\,y\right)\dd y_1\dd y_2\right),\]
\[R_2 = \frac{1}{\sqrt{\det(\Sigma^{\diag})}}\iint_{\R^2}|y_1||y_2|\left[\exp\left(-\frac{1}{2}\,y^T(\Sigma_c)^{-1}\,y\right)-\exp\left(-\frac{1}{2}\,y^T( \Sigma^{\diag})_c)^{-1}\,y\right)\right]\dd y_1\dd y_2.\]
Estimate \eqref{eqbis:05} directly implies
\[R_1 = O\left(\|\Sigma_{\off}(s,t)\|^2\right).\numberthis\label{eqbis:08}\]
A Taylor expansion of the exponential function and relations \eqref{eqbis:06} and \eqref{eqbis:07} yields
\[ R_2 = \frac{2h}{\sqrt{\det(\Sigma^{\diag}}}\iint_{\R^2}|y_1||y_2|y_1y_2\exp\left(-\frac{1}{2}\,y^T( \Sigma^{\diag})_c)^{-1}\,y\right)\dd y_1\dd y_2 + O(\|\Sigma^{\off}\|^2).\]
By parity the double integral is zero, thus
\[R_2=O\left(\|\Sigma_{\off}(s,t)\|^2\right).\numberthis\label{eqbis:09}\]
The Lagrange rest theorem implies that all the $O$ appearing above are uniform for $\Sigma$ in the compact set $K$. Gathering estimates \eqref{eqbis:08} and \eqref{eqbis:09}, we deduce the existence of a constant $C_K$ such that
\[\forall \Sigma\in K,\quad\left|\tilde{\rho}_2(\Sigma) - \tilde{\rho}_2(\Sigma^{\diag})\right|\leq C_K\|\Sigma^{\off}\|^2.\]
\end{proof}
\begin{corollary}
\label{lemmabis8}
There are positive constants $\varepsilon$ and $C$ depending only on the constant $M$ defined in \eqref{eqbis:02}, such that for all $s,t\in I$ satisfying $\max_{u,v\in\{0,1\}}|r^{(u,v)}(s,t)|\leq \varepsilon$,
\[\det \Sigma(s,t) \geq \frac{1}{2C^2}\quad\quand\quad |\rho_2(s,t)-\rho_1(s)\rho_1(t)|\leq C\left(\sup_{u,v\in\{0,1\}}|r^{(u,v)}(s,t)|\right)^2.\]
\end{corollary}
\begin{proof}
Given the equality \eqref{eqbis:03}, the upper bound assumption \eqref{eqbis:02} and the regularity of the determinant, we deduce the existence of a constant $C$ depending only on the constant $M$, such that
\[|\det \Sigma(s,t)-\det\Omega(s)\det \Omega(t)|\leq C\|\Sigma_{\off}(s,t)\| = C\!\!\!\max_{u,v\in\{0,1\}}|r^{(u,v)}(s,t)|.\]
The lower bound assumption on the determinant \eqref{eqbis:02} implies the existence of a positive number $\varepsilon$ depending only on $M$ such that for all $s,t\in I$ with $\max_{u,v\in\{0,1\}}|r^{(u,v)}(s,t)|\leq \varepsilon$ one has
\[\det \Sigma(s,t) \geq \frac{1}{2M^2}.\numberthis\label{eqbis:10}\]
From now we fix two such real numbers $s$ and $t$. The above lower bound \eqref{eqbis:10} implies that the matrix $\Sigma(s,t)$ lives in a compact subset $K_M$ of the space of symmetric positive definite matrices, that depends only on the constant $M$.  The previous Lemma \ref{lemmabis4} applied to the matrix $\Sigma(s,t)$ and the relations \eqref{eqbis:11} imply the existence of a constant $C$ depending only on $M$ such that
\[|\rho_2(s,t)-\rho_1(s)\rho_1(t)|\leq C\left(\max_{u,v\in\{0,1\}}|r^{(u,v)}(s,t)|\right)^2.\]
\end{proof}

\section{Proof of the main theorems}\label{secbis3}
\subsection{Asymptotics of the Kac density}
Let $(X_n)_n$ be a sequence of processes on $I$ satisfying the hypotheses of Theorem \ref{theorembis2}. We recall the definition for $n\in\N$ of the process $Y_n = X_n(\,.\,/n)$ with covariance function $r_n$, and $Y_\infty$ is a centered stationary Gaussian process with covariance function $r_\infty$. We define for $n\in\N\cup\{+\infty\}$
\[\Omega_n(s) := \Cov(Y_n(s),Y'_n(s))\quand \Sigma_n(s,t) = \Cov(Y_n(s),Y_n(t),Y'_n(s),Y'_n(t)).\]
By stationarity, $\Omega_\infty$ is constant and the matrix $\Sigma_n(s,t)$ depends only on the difference $s-t$. The following Lemma \ref{lemmabis3} shows that the process $Y_n$ satisfies the bound \eqref{eqbis:02}, uniformly for $n$ large enough.
\begin{lemma}
\label{lemmabis3}
There is a constant $C$ such that for $s,t\in nI$ and $u,v\in\{0,1,2\}$ one has
\[|r_n^{(u,v)}(s,t)|\leq C.\]
For all parameter $\eta>0$, there is a positive constant $C_\eta$ and a rank $n_0$ such that for $n\geq n_0$ and for all $s,t\in nI$ with $\dist(s,t)>\eta$,
\[\det \Sigma_n(s,t) \geq C_\eta.\]
\end{lemma}
\begin{proof}
The uniform convergence of hypothesis $(H_1)$ on the process $X_n$ ensures that the quantities $\E[Y_n(s)^2]$, $\E[Y_n'(s)^2]$ and $\E[Y''_n(s)^2]$ are bounded by some constant $C$. By Cauchy-Schwartz inequality, for $u,v\in\{0,1,2\}$ and $s,t\in nI$ one has $|r_n^{(u,v)}(s,t)|\leq C$. For $s\in \R$, the random variables $Y_\infty(s)$ and $Y'_\infty(s)$ are decorrelated by stationarity, thus
\[\det\Omega_\infty = -r_\infty(0)r''_\infty(0) >0.\]
We deduce the following uniform convergence on $x\in I$
\[\lim_{n\rightarrow+\infty} \Omega_n(nx) = \psi(x)\Omega_\infty.\numberthis\label{eqbis:12}\]
The function $\psi$ is bounded from below by a positive constant $C_\psi$, and the convergence \eqref{eqbis:12} implies the existence of rank $n_0$ independent of $x$ such that
\[\forall n\geq n_0,\forall x\in I,\quad\det[\Omega_n(nx)]\geq \frac{C_\psi^2 \det\Omega_\infty}{2}>0.\numberthis\label{eqbis:24}\]
In particular, the function $Y_n$ satisfies the bounds \eqref{eqbis:02} for some constant $M$ independent of $n$. Now the function $g$ of hypothesis $(H_2)$ in Theorem \ref{theorembis2} decreases to zero at infinity. Given $\varepsilon>0$, there is a constant $T_\varepsilon$ such that for $s,t$ satisfying $\dist(s,t)>T_\varepsilon$, one has
\[\max_{u,v\in\{0,1\}}|r_n^{(u,v)}(s,t)|\leq \varepsilon.\]
One can then apply Corollary \ref{lemmabis8}, to deduce the existence of a constant $T$ independent of $n$, such that for all $s,t$ satisfying $\dist(s,t)>T$, it holds that 
\[\det \Sigma_n(s,t) \geq \frac{1}{2M^2}.\numberthis\label{eqbis:15}\]
\jump
Let $\eta>0$. Since the process $Y_\infty$ is stationary and the support of its spectral measure has an accumulation point, then (see \cite[Ex. ~3.5]{Aza09}) the covariance matrix $\Sigma_\infty(s,t)$ is nondegenerate for $s\neq t$. By compactness, one can find a positive constant $C_\eta$ such that for all $s,t\in nI$ with $\eta\leq \dist(s,t)\leq T$, one has $\det\Sigma_\infty(s-t)\geq C_\eta$.
The uniform convergence of $r^{(u,v)}_n$ towards $r^{(u,v)}_\infty$ and then implies that for $n$ greater than some rank $n_0$, 
\[\det\Sigma_n(s,t)\geq \frac{C_\eta}{2}.\numberthis\label{eqbis:13}\]
Gathering \eqref{eqbis:15} and \eqref{eqbis:13} we deduce the second assertion.
\end{proof}

For $n\in \N\cup\{+\infty\}$ we define $\rho_{1,n}$ and $\rho_{2,n}$ the Kac densities associated to the process $Y_n$. By stationarity, the function $\rho_{1,\infty}$ is a constant and the function $\rho_{2,\infty}$ depends only on the difference $s-t$. We deduce the following corollary.
\begin{corollary}
\label{lemmabis5}
Let $\eta>0$. We have the following uniform convergences, uniformly in $x\in I$ and $s,t$ in a compact set of $\R$ with $\dist(s,t)>\eta$.
\[\lim_{n\rightarrow+\infty} \rho_{1,n}(nx) = \rho_{1,\infty}\quand\lim_{n\rightarrow+\infty} \rho_{2,n}(nx+s,nx+t) = \rho_{2,\infty}(s-t).\]
\end{corollary}
\begin{proof}
Hypothesis $(H_1)$ and the previous Lemma \ref{lemmabis3} implies the following uniform convergences uniformly in $x\in I$ and $s,t$ in a compact set of $\R$ with $|s-t|>\eta$
\[\lim_{n\rightarrow+\infty} \Omega_n(nx) = \psi(x)\Omega_\infty\quand\lim_{n\rightarrow+\infty} \Sigma_n(nx+s,nx+t) = \psi(x)\Sigma_\infty(s-t).\numberthis\label{eqbis:16}\]
The functions $\tilde{\rho}_1$ (resp. $\tilde{\rho}_2$) is continuous on the space of symmetric positive definite matrices, which combined with the convergence \eqref{eqbis:16} directly implies the conclusion of the corollary. Note that by a change of variable, the limit does not depend on the function $\psi$. 
\end{proof}

The following Lemma \ref{lemmabis6} establishes a decay property for the Kac density, whose rate depends on the function $g$ of hypothesis $(H_2)$ of Theorem \ref{theorembis2}.
\begin{lemma}
\label{lemmabis6}
There is a constant $C$ and a rank $n_0$ such that for $n\geq n_0$ and $s,t\in nI$,
\[\left|\rho_{2,n}(s,t) - \rho_{1,n}(s)\rho_{1,n}(t)\right| \leq Cg^2(\dist(s,t)).\]
\end{lemma}
\begin{proof}
According to Corollary \ref{lemmabis8} and the proof of Lemma \ref{lemmabis3}, there are constants $T$ and $C$ independent of $n$ such that when $\dist(s,t)\geq T$ it holds that
\[\left|\rho_{2,n}(s,t) - \rho_{1,n}(s)\rho_{1,n}(t)\right|\leq Cg^2(\dist(s,t)).\] 
When $\dist(s,t)\leq T$, the functions $\rho_{1,n}$ and $\rho_{2,n}$ are bounded by a constant independent of $n$, according to Lemmas \ref{lemmabis9} and \ref{lemmabis3}. Since the function $g$ is assumed to be continuous and positive, it is bounded below by a positive constant on compact subsets of $\R$ and the conclusion follows.
\end{proof}

\subsection{Proof of Theorem \ref{theorembis2}}
We identify the interval $I$ with $]a,b[$ and we write
\begin{align*}
\Var(Z_{X_n}(I)) &= \Var(Z_{Y_n}(nI))\\
&=n\int_{a}^{b}\!\rho_{1,n}(nx)\dd x + n\int_{a}^{b}\!\int_{n(a-x)}^{n(b-x)} \left(\rho_{2,n}(nx,nx+\tau) - \rho_{1,n}(nx)\rho_{1,n}(nx+\tau)\right)\dd \tau\dd x.
\end{align*}
For the first term, Corollary \ref{lemmabis5} implies that
\[\lim_{n\rightarrow+\infty} \int_{a}^{b}\rho_{1,n}(nx)\dd x = |b-a|\rho_{1,\infty}.\]
For the second term, we write
\begin{align*}
R_n = &\int_{a}^{b}\int_{n(a-x)}^{n(b-x)} \left(\rho_{2,n}(nx,nx+\tau) - \rho_{1,n}(nx)\rho_{1,n}(nx+\tau)\right)\dd \tau\dd x\\
&\quad-\quad|b-a|\int_\R \left(\rho_{2,\infty}(\tau) - \rho_{1,\infty}^2\right)\dd \tau.
\end{align*}
We fix two positive constants $\eta>0$ and $A$. We split $R_n$ into four parts, as
\[R_n = R_{1,n}^{A,\eta}- R_{2,n}^{A,\eta}+ R_{3,n}^{A,\eta} - R_{3,\infty}^{A,\eta},\]
with
\begin{align*}
R_{1,n}^{A,\eta} &:= \int_{a}^{b}\int_{n(a-x)}^{n(b-x)}\one_{\eta\leq |\tau|\leq A} \left[\left(\rho_{1,n}(nx)\rho_{1,n}(nx+\tau)-\rho_{1,\infty}^2\right)\right]\dd\tau\dd x,\\
R_{2,n}^{A,\eta} &:= \int_{a}^{b}\int_{n(a-x)}^{n(b-x)}\one_{\eta\leq |\tau|\leq A} \left[\left(\rho_{2,n}(nx,nx+\tau)-\rho_{2,\infty}(\tau)\right)\right]\dd \tau\dd x,\\ 
R_{3,n}^{A,\eta} &:= \int_{a}^{b}\int_{n(a-x)}^{n(b-x)}\one_{\{|\tau|\geq A\}\cup\{|\tau|\leq \eta\}} \left(\rho_{2,n}(nx,nx+\tau) - \rho_{1,n}(nx)\rho_{1,n}(nx+\tau)\right)\dd \tau\dd x,\\
R_{3,\infty}^{A,\eta} &:= |b-a|\int_\R \one_{\{|\tau|\geq A\}\cup\{|\tau|\leq \eta\}}(\rho_{2,\infty}(\tau) - \rho_{1,\infty}^2)\dd \tau.
\end{align*}
Corollary \ref{lemmabis5} directly implies that $R_{1,n}^{A,\eta}$ and $R_{2,n}^{A,\eta}$ converge towards $0$ when $n$ goes to infinity. The bound given by Lemma \ref{lemmabis6} implies that for some constant $C$ independent of $A$ and $\eta$,
\begin{align*}
|R_{3,n}^{A,\eta}| + |R_{3,\infty}^{A,\eta}| \leq C|b-a|\int_\R \one_{\{|\tau|\geq A\}\cup\{|\tau|\leq \eta\}}g^2(\tau)\dd \tau.\numberthis\label{eqbis:32}
\end{align*}
Gathering the estimates for $R_{3,n}^{A,\eta}$, $R_{3,n}^{A,\eta}$, $R_{3,n}^{A,\eta}$ and $R_{3,\infty}^{A,\eta}$, we deduce that for some constant $C$,
\[\limsup_{n\rightarrow+\infty} |R_n| \leq C|b-a|\int_\R \one_{\{|\tau|\geq A\}\cup\{|\tau|\leq \eta\}} g^2(\tau)\dd \tau.\]
The function $g$ is assumed to be square integrable. Letting $A$ go to infinity and $\eta$ go to zero, we deduce that $\lim_{n\rightarrow+\infty} R_n = 0$, from which follows the following convergence
\[\lim_{n\rightarrow +\infty} \frac{\Var(Z_{X_n}(I))}{n} = \mathrm{length}(I)\,\gamma_2,\]
where
\[\gamma_2 := \int_\R \left(\rho_{2,\infty}(\tau) - \rho_{1,\infty}^2\right)\dd \tau + \rho_{1,\infty}.\]
It remains to show the positivity of the constant $C_\infty$. The above proof shows in fact that 
\[\gamma_2 = \lim_{n\rightarrow +\infty} \frac{\Var(Z_{Y_\infty}[0,n])}{n},\]
and it has been shown (see \cite{Anc20}) that $C_\infty>0$ for a large class of processes including $Y_\infty$.

\subsection{Proof of Theorem \ref{theorembis1}}

In the following, we consider the sequence of trigonometric Gaussian polynomials $(X_n)_{n}$ defined in \eqref{eqbis:17}. Assuming the hypotheses of Theorem \ref{theorembis1}, we show that this sequence of processes satisfies hypotheses $(H_1)$ and $(H_2)$ of Theorem \ref{theorembis2}, from which follows the conclusion of Theorem \ref{theorembis1}. Following \cite{Ang21}, the next computation gives an integral expression for the covariance function $r_n$ of $Y_n = X_n(./n)$.
\begin{align*}
r_n(s,t)  :\!&= \E[Y_n(s)Y_n(t)]\\
&= \frac{1}{n}\sum_{k,l=1}^n \rho(k-l)\cos\left(\frac{ks-lt}{n}\right)\\
&=\frac{1}{2\pi n}\int_0^{2\pi}\mathrm{Re}\left(\sum_{k,l=1}^ne^{-i(k-l)y}e^{\frac{iks-ilt}{n}}\right)\dd\mu(y)\\
&= \cos\left(\frac{n+1}{2n}(s-t)\right)\frac{1}{2\pi}\int_0^{2\pi} K_n\left(\frac{s}{n}-y,\frac{t}{n}-y\right)\psi(y)\dd y,\numberthis\label{eqbis:36}
\end{align*}
where $K_n$ is the two points Fejér kernel
\[K_n(x,y) = \frac{1}{n}\frac{\sin\left(\frac{nx}{2}\right)}{\sin\left(\frac{x}{2}\right)}\frac{\sin\left(\frac{ny}{2}\right)}{\sin\left(\frac{y}{2}\right)}.\]
In the case where $\rho(k-l) = \delta_{k,l}$, the measure $\mu$ is the normalized Lebesgue measure on $[-\pi,\pi]$. In that case, we denote by $r_{0,n}$ its covariance function, whose expression is given by
\[r_{0,n}(s,t) = \frac{1}{2n}\left[\frac{\sin\left(\left(\frac{2n+1}{2n}\right)(s-t)\right)}{\sin\left(\frac{s-t}{2n}\right)}-1\right].\numberthis\label{eqbis:01}\]

By assumption the spectral measure $\mu$ has a continuous and positive density $\psi$ on $\T$. The following two lemmas show that the covariance function $r_n$ satisfies hypotheses $(H_1)$ and $(H_2)$ of Theorem \ref{theorembis2} with function $g(t) = C/(1+|t|^\alpha)$, for $1/2<\alpha<1$, and $Y_\infty$ a stationary Gaussian process with $\sinc$ covariance function.\jump

\begin{lemma}
\label{lemmabis2}
Let $u,v\geq 0$. Uniformly for $x\in \T$ and $s,t$ in compact subsets of $\R$,
\[\lim_{n\rightarrow+\infty} r_n^{(u,v)}(nx+s,nx+t) = \psi(x)(-1)^v\sinc^{(u+v)}(s-t).\]
\end{lemma}

\begin{proof}
Let us first remark that the covariance function $r_n$ is a trigonometric polynomial and can thus be extended to an analytic function on $\C$. We will prove that the conclusion of Lemma \ref{lemmabis2} holds when $s$ and $t$ belong to a compact subset of $\C$. By analyticity, it suffices then to prove the lemma for $u=v=0$. We have
\begin{align*}
r_n(nx + s,nx+t) &= I_n^x(s,t) + \psi(x)r_{0,n}(nx + s,nx+t),
\end{align*}
where 
\[I_n^x(s,t)  = \cos\left(\frac{n+1}{2n}(s-t)\right)\frac{1}{2\pi}\int_{-\pi}^{\pi} K_n\left(\frac{s}{n}-y,\frac{t}{n}-y\right)\left[\psi\left(x+y\right)-\psi(x)\right]\dd y.\]
From expression \eqref{eqbis:01}, one has uniformly for $x\in\T$ and $s,t$ in compact subsets of $\C$,
\[r_{0,n}(nx+s,nx+t) = \sinc(s-t) + O\left(\frac{1}{n}\right).\]
It remains to prove that the quantity $I_n^x(s,t)$ converges towards $0$ uniformly on $x\in\T$ and $s,t$ in compact sets of $\C$. Let $A>1$ and $K = B(0,A-1)$ the disk centered in $0$ of radius $A-1$ in $\C$. Denoting by $\omega_\psi$ the uniform modulus of continuity of the spectral density $\psi$, we have 
\begin{align*}
|I_n^x(s,t)| 
&\leq \frac{1}{2\pi n^2}\int_{-n\pi}^{n\pi} \left|\frac{\sin\left(\frac{s-y}{2}\right)}{\sin\left(\frac{s-y}{2n}\right)}\frac{\sin\left(\frac{t-y}{2}\right)}{\sin\left(\frac{t-y}{2n}\right)}\right|\omega_\psi\left(\frac{y}{n}\right)\dd y\\
&\leq R_1 + R_2,
\end{align*}
where 
\[R_1 = \frac{1}{2\pi n^2}\int_{-A}^A \left|\frac{\sin\left(\frac{s-y}{2}\right)}{\sin\left(\frac{s-y}{2n}\right)}\frac{\sin\left(\frac{t-y}{2}\right)}{\sin\left(\frac{t-y}{2n}\right)}\right|\omega_\psi\left(\frac{y}{n}\right)\dd y,\]
and
\[R_2 = \frac{1}{2\pi n^2}\int_{-n\pi}^{n\pi} \one_{\lbrace|y|\geq A\rbrace}\left|\frac{\sin\left(\frac{s-y}{2}\right)}{\sin\left(\frac{s-y}{2n}\right)}\frac{\sin\left(\frac{t-y}{2}\right)}{\sin\left(\frac{t-y}{2n}\right)}\right|\omega_\psi\left(\frac{y}{n}\right)\dd y.\]
The term $R_1$ is bounded by
\[R_1 \leq \frac{1}{2\pi}\int_{-A}^A \omega_\psi\left(\frac{y}{n}\right)\dd y.\]
Since the spectral density is (uniformly) continuous on $\T$, the quantity $R_1$ converges towards zero as $n$ goes to infinity, uniformly on $x\in\T$ and $s,t\in K$. For the term $R_2$ we use the following inequalities, valid for $\mathrm{Re}(z)\in[-5\pi/6,5\pi/6]$:
\[\frac{3}{5\pi}|\mathrm{Re}(z)|\leq |\sin(\mathrm{Re}(z))|\leq |\sin(z)|\numberthis\label{eqbis:22}.\]
There is a rank $n_0$ depending only on the compact subset $K$ such that, for all $n\geq n_0$, $s\in K$ and $y\in [-n\pi,n\pi]$,
\[-\frac{5\pi}{6}\leq \frac{\mathrm{Re}(s)-y}{2n}\leq \frac{5\pi}{6}\quand -\frac{5\pi}{6}\leq \frac{\mathrm{Re}(t)-y}{2n}\leq \frac{5\pi}{6}.\]
It follows from the series of inequalities \eqref{eqbis:22} that there is a constant $C$ such that for $n\geq n_0$,
\begin{align*}
R_2 &\leq \frac{\cosh(A)^2}{2\pi n^2}\int_{-n\pi}^{n\pi} \one_{\lbrace|y|\geq A\rbrace}\frac{1}{\left|\sin\left(\frac{s-y}{2n}\right)\right|}\frac{1}{\left|\sin\left(\frac{t-y}{2n}\right)\right|}\omega_\psi\left(\frac{y}{n}\right)\dd y\\
&\leq C\cosh(A)^2\int_{-n\pi}^{n\pi} \one_{\lbrace|y|\geq A\rbrace}\frac{1}{|\mathrm{Re}(s)-y||\mathrm{Re}(t)-y|}\omega_\psi\left(\frac{y}{n}\right)\dd y\\
&\leq C\cosh(A)^2\int_{-\infty}^\infty \one_{\lbrace|y|\geq A\rbrace}\frac{1}{\left(|y|-A+1\right)^2}\,\omega_\psi\left(\frac{y}{n}\right)\dd y.
\end{align*}
The quantity $R_2$ thus converges towards $0$ as $n$ goes to infinity, uniformly on $x\in\T$ and $s,t\in K$.
\end{proof}

\begin{lemma}
\label{lemmabis1}
Let $u,v\geq 0$ and $0<\alpha<1$. There is a constant $C$ such that
\begin{align*}
\forall s,t\in n\T,\quad |r_n^{(u,v)}(s,t)|\leq \frac{C}{1+\dist(s,t)^\alpha}.
\end{align*}
\end{lemma}
\begin{proof}
Let $s,t\in n\T$. According to the previous Lemma \ref{lemmabis2}, the function $r_n$ and its derivatives are uniformly bounded, we can thus assume that $\dist(s,t)\geq 4$. By Cauchy integral formula, there is a constant $C$ such that
\begin{align*}
|r_n^{(a,b)}(s,t)| \leq C\sup_{|w|\leq 1}\sup_{|z|\leq 1} \left|r_n\left(s+w,t+z\right)\right|.
\end{align*}
Let $w,z$ be complex numbers such that $|w|\leq 1$ and $|z|\leq 1$. Using the explicit formula \eqref{eqbis:36} for $r_n$ we obtain
\begin{align*}
r_n(s+w,t+z) &= \frac{1}{2\pi n}\cos\left(\frac{n+1}{2n}(s-t+w-z)\right)\!\!\int_{-\pi}^\pi\! \frac{\sin\left(\frac{ny+s+w}{2}\right)}{\sin\left(\frac{y+\frac{s+w}{n}}{2}\right)}\frac{\sin\left(\frac{ny+t+z}{2}\right)}{\sin\left(\frac{y+\frac{t+z}{n}}{2}\right)}\psi(y)\dd y.
\end{align*}
Using the fact that the cosine function is bounded on a horizontal complex strip,
\begin{align*}
|r_n(s+w,t+z)| &\leq \frac{C\|\psi\|_\infty}{2\pi n}\int_{-\pi}^\pi \left|\frac{\sin\left(\frac{ny+s+w}{2}\right)}{\sin\left(\frac{y+\frac{s+w}{n}}{2}\right)}\frac{\sin\left(\frac{ny+t+z}{2}\right)}{\sin\left(\frac{y+\frac{t+z}{n}}{2}\right)}\right|\dd y.
\end{align*}
Let $\delta = \dist(s,t)/2$. Up to translating $s$ and $t$ by $\pm 2\pi n$ and exchanging $s$ and $t$, we can assume that $\delta = \frac{t-s}{2}$. We then make the change of variable $x = y + \frac{t+s}{2n}$ to obtain
\[\int_{-\pi}^\pi \left|\frac{\sin\left(\frac{ny+s+w}{2}\right)}{\sin\left(\frac{y+\frac{s+w}{n}}{2}\right)}\frac{\sin\left(\frac{ny+t+z}{2}\right)}{\sin\left(\frac{y+\frac{t+z}{n}}{2}\right)}\right|\dd y = 
\int_{-\pi}^\pi\left|\frac{\sin\left(\frac{nx+\delta+w}{2}\right)}{\sin\left(\frac{x+\frac{\delta+w}{n}}{2}\right)}\frac{\sin\left(\frac{nx-\delta+z}{2}\right)}{\sin\left(\frac{x+\frac{-\delta+z}{n}}{2}\right)}\right|\dd x.\]
This last integral splits into two integrals $I_1$ and $I_2$ defined by 
\[I_1 := \int_{0}^\pi \left|\frac{\sin\left(\frac{nx+\delta+w}{2}\right)}{\sin\left(\frac{x+\frac{\delta+w}{n}}{2}\right)}\frac{\sin\left(\frac{nx-\delta+z}{2}\right)}{\sin\left(\frac{x+\frac{-\delta+z}{n}}{2}\right)}\right|\dd x\quand I_2 :=\int_{-\pi}^0 \left|\frac{\sin\left(\frac{nx+\delta+w}{2}\right)}{\sin\left(\frac{x+\frac{\delta+w}{n}}{2}\right)}\frac{\sin\left(\frac{nx-\delta+z}{2}\right)}{\sin\left(\frac{x+\frac{-\delta+z}{n}}{2}\right)}\right|\dd x.\]
Both term can be treated the exact same way. We have by Hölder inequality with $0<\alpha<1$,
\begin{align*}
I_1&\leq \left(\int_0^\pi \left|\frac{\sin\left(\frac{nx+\delta+w}{2}\right)}{\sin\left(\frac{x+\frac{\delta+w}{n}}{2}\right)}\right|^\frac{1}{1-\alpha}\dd x\right)^{1-\alpha}
\left(\int_0^\pi \left|\frac{\sin\left(\frac{nx-\delta+z}{2}\right)}{\sin\left(\frac{x+\frac{-\delta+z}{n}}{2}\right)}\right|^\frac{1}{\alpha}\dd x\right)^\alpha.\numberthis\label{eqbis:38}
\end{align*}
For the left integral in \eqref{eqbis:38}, we make use of the following inequalities, which are consequences of inequalities \eqref{eqbis:22}, and the fact that $|w|\leq 1$ and $\delta\geq 2$
\[\left|\sin\left(\frac{x+\frac{\delta+w}{n}}{2}\right)\right|\geq \frac{3}{10\pi}\left(x + \frac{\delta+\mathrm{Re}(u)}{n}\right)\geq\frac{3}{10\pi}\left(x+\frac{\delta}{2n}\right),\]
to get
\begin{align*}
\left(\int_0^\pi \left|\frac{\sin\left(\frac{nx+\delta+w}{2}\right)}{\sin\left(\frac{x+\frac{\delta+w}{n}}{2}\right)}\right|^\frac{1}{1-\alpha}\dd x\right)^{1-\alpha}&\leq
C\left(\int_0^\infty \frac{\dd x}{\left(x+\frac{\delta}{2n}\right)^\frac{1}{1-\alpha}}\right)^{1-\alpha}\\
&\leq C\left(\frac{n}{\delta}\right)^\alpha.\numberthis\label{eqbis:37}
\end{align*}
For the right integral in \eqref{eqbis:38}, we make the change of variable $y=nx-\delta+\mathrm{Re}(z)$ and we use the inequality
\[|x+iy|\leq C|\sin(x+iy)|,\]
valid for $x\in [-5\pi/6,5\pi/6]$ and $y\in [-1/4,1/4]$, to get
\begin{align*}
\left(\int_0^\pi \left|\frac{\sin\left(\frac{nx-\delta+z}{2}\right)}{\sin\left(\frac{x+\frac{-\delta+z}{n}}{2}\right)}\right|^\frac{1}{\alpha}\dd y\right)^\alpha&\leq
n^{-\alpha}\left(\int_{\mathrm{Re}(z)-\delta}^{n\pi + \mathrm{Re}(z)-\delta} \left|\frac{\sin\left(\frac{y+i\,\mathrm{Im}(v)}{2}\right)}{\sin\left(\frac{y+i\,\mathrm{Im}(v)}{2n}\right)}\right|^\frac{1}{\alpha}\dd y\right)^\alpha\\
&\leq 2Cn^{1-\alpha}\left(\int_{-\infty}^\infty \left|\frac{\sin\left(\frac{y+i\,\mathrm{Im}(v)}{2}\right)}{y+i\,\mathrm{Im}(v)}\right|^\frac{1}{\alpha}\dd y\right)^\alpha\\
&\leq C'n^{1-\alpha},\numberthis\label{eqbis:39}
\end{align*}
where in the last inequality we used the fact that the integrand is uniformly bounded in a neighborhood of the origin, and that $\frac{1}{\alpha}>1$ so the integrand is also integrable near $\pm\infty$. Plugging estimates \eqref{eqbis:37} and \eqref{eqbis:39} into inequality \eqref{eqbis:38} we obtain for some constant $C$ that
\[|r_n^{(u,v)}(nx+s,nx+t)
)|\leq \frac{C\|\psi\|_\infty}{2\pi n}(I_1+I_2) \leq \frac{2C\|\psi\|_\infty}{2\pi n}\left(\frac{n}{\delta^\alpha}\right) \leq  \frac{C}{\dist(s,t)^\alpha}.\]
\end{proof}

\printbibliography

@misc{Anc20,
      title={Zeros of smooth stationary Gaussian processes}, 
      author={Ancona, Michele and Letendre, Thomas},
      year={2020},
      eprint={2007.03240},
      archivePrefix={arXiv},
      primaryClass={math.PR}
}

@article {Ang19,
    AUTHOR = {Angst, J\"{u}rgen and Dalmao, Federico and Poly, Guillaume},
     TITLE = {On the real zeros of random trigonometric polynomials with
              dependent coefficients},
   JOURNAL = {Proc. Amer. Math. Soc.},
  FJOURNAL = {Proceedings of the American Mathematical Society},
    VOLUME = {147},
      YEAR = {2019},
    NUMBER = {1},
     PAGES = {205--214},
      ISSN = {0002-9939},
   MRCLASS = {42A05 (30C15)},
  MRNUMBER = {3876743},
       DOI = {10.1090/proc/14216},
       URL = {https://doi.org/10.1090/proc/14216},
}

@article {Ang21,
    AUTHOR = {Angst, J\"{u}rgen and Pautrel, Thibault and Poly, Guillaume},
     TITLE = {Real zeros of random trigonometric polynomials with dependent coefficients},
       URL = {https://arxiv.org/abs/2102.09653},
      YEAR = {2021},
}

@book {Aza09,
    AUTHOR = {Aza\"{i}s, Jean-Marc and Wschebor, Mario},
     TITLE = {Level sets and extrema of random processes and fields},
 PUBLISHER = {John Wiley \& Sons, Inc., Hoboken, NJ},
      YEAR = {2009},
     PAGES = {xii+393},
      ISBN = {978-0-470-40933-6},
   MRCLASS = {60-02 (60E15 60G05 60G15 60G60 60G70)},
  MRNUMBER = {2478201},
MRREVIEWER = {Anna Amirdjanova},
       DOI = {10.1002/9780470434642},
       URL = {https://doi.org/10.1002/9780470434642},
}

@article {Aza16,
    AUTHOR = {Aza\"{i}s, Jean-Marc and Dalmao, Federico and Le\'{o}n, Jos\'{e} R.},
     TITLE = {C{LT} for the zeros of classical random trigonometric
              polynomials},
   JOURNAL = {Ann. Inst. Henri Poincar\'{e} Probab. Stat.},
  FJOURNAL = {Annales de l'Institut Henri Poincar\'{e} Probabilit\'{e}s et
              Statistiques},
    VOLUME = {52},
      YEAR = {2016},
    NUMBER = {2},
     PAGES = {804--820},
      ISSN = {0246-0203},
   MRCLASS = {60G15 (42C10 60F05)},
  MRNUMBER = {3498010},
MRREVIEWER = {Maxim L. Yattselev},
       DOI = {10.1214/14-AIHP653},
       URL = {https://doi.org/10.1214/14-AIHP653},
}

@article {Bal19,
    AUTHOR = {Bally, Vlad and Caramellino, Lucia and Poly, Guillaume},
     TITLE = {Non universality for the variance of the number of real roots
              of random trigonometric polynomials},
   JOURNAL = {Probab. Theory Related Fields},
  FJOURNAL = {Probability Theory and Related Fields},
    VOLUME = {174},
      YEAR = {2019},
    NUMBER = {3-4},
     PAGES = {887--927},
      ISSN = {0178-8051},
   MRCLASS = {60G50 (60F05)},
  MRNUMBER = {3980307},
       DOI = {10.1007/s00440-018-0869-2},
       URL = {https://doi.org/10.1007/s00440-018-0869-2},
}

@article {Gra11,
    AUTHOR = {Granville, Andrew and Wigman, Igor},
     TITLE = {The distribution of the zeros of random trigonometric
              polynomials},
   JOURNAL = {Amer. J. Math.},
  FJOURNAL = {American Journal of Mathematics},
    VOLUME = {133},
      YEAR = {2011},
    NUMBER = {2},
     PAGES = {295--357},
      ISSN = {0002-9327},
   MRCLASS = {60G99 (60F05)},
  MRNUMBER = {2797349},
       DOI = {10.1353/ajm.2011.0015},
       URL = {https://doi.org/10.1353/ajm.2011.0015},
}

@article {Lub21,
    AUTHOR = {Lubinsky, Doron S. and Pritsker, Igor E.},
     TITLE = {Variance of real zeros of random orthogonal polynomials},
   JOURNAL = {J. Math. Anal. Appl.},
  FJOURNAL = {Journal of Mathematical Analysis and Applications},
    VOLUME = {498},
      YEAR = {2021},
    NUMBER = {1},
     PAGES = {124954},
      ISSN = {0022-247X},
   MRCLASS = {30C15 (60G99)},
  MRNUMBER = {4202193},
       DOI = {10.1016/j.jmaa.2021.124954},
       URL = {https://doi.org/10.1016/j.jmaa.2021.124954},
}
\end{document}